\documentstyle[12pt,amsfonts,amssymb,leqno]{article}

\newfont\got{eufm10}

\newtheorem{proposition}{Proposition}[section]
\newtheorem{thm}[proposition]{Theorem}

\newtheorem{lemma}[proposition]{Lemma}
\newtheorem{defn}[proposition]{Definition}

\renewcommand{\thefootnote}{\alph{footnote}}

\newcounter{secnum}

\newcommand{\forces}{\ \ |\!\!\!|\!- \ }

\begin{document}
\setcounter{section}{+0}

\begin{center}
{\Large \bf Bounded Martin's Maximum is stronger}
\end{center}
\begin{center}
{\Large \bf than the Bounded Semi-proper Forcing Axiom}
\end{center}

\begin{center}
{\large Ralf Schindler}
\renewcommand{\thefootnote}{arabic{footnote}}
\end{center}
\begin{center} 
{\footnotesize
{\it Institut f\"ur Formale Logik, Universit\"at Wien, 1090 Wien, Austria}} 
\end{center}

\begin{center}
{\tt rds@logic.univie.ac.at}

{\tt http://www.logic.univie.ac.at/${}^\sim$rds/}\\
\end{center}

\begin{abstract}
\noindent
We show that if Bounded Martin's Maximum (${\sf BMM}$)
holds then for every
$X \in V$ there is an inner model with a strong cardinal containing
$X$. In particular, by \cite{gosh}, ${\sf BMM}$ is strictly stronger
consistency-wise than the Bounded Semi-Proper Forcing Axiom
(${\sf BSPFA}$).
\end{abstract}

\section{Introduction.}

Shelah has shown that the Semi-Proper Forcing Axiom (${\sf SPFA}$) is
equivalent with Martin's Maximum (${\sf MM}$). It was an open problem 
to decide whether the same holds true 
at least consistency-wise
for the bounded versions
of these axioms, i.e., 
to decide whether the Bounded Semi-Proper Forcing Axiom (${\sf BSPFA}$)
is really or only apparently weaker than Bounded 
Martin's Maximum (${\sf BMM}$). In this paper we shall solve this problem
by showing that ${\sf BMM}$ yields the existence of inner models with strong
cardinals; in fact, we shall prove:

\begin{thm}\label{main}
Suppose that ${\sf BMM}$ holds. Then for every
$X \in V$ there is an inner model with a strong cardinal containing
$X$. 
\end{thm}

The key technical lemma which will give Theorem \ref{main} is
Lemma \ref{key}; this lemma is shown by designing a refined
$K$-version of Jensen's ``reshaping'' (the paper \cite{K}
contains such a version which is almost good enough for the present
purpose).\footnote{The author would like to thank David Asper\'o
for a pivotal discussion about ${\sf BMM}$.}

By \cite{gosh}, ${\sf BSPFA}$ is equiconsistent with a reflecting cardinal, 
which lives consistency-wise between inaccessible and Mahlo cardinals.
Theorem \ref{main} therefore implies that ${\sf BMM}$ is consistency-wise
strictly stronger than ${\sf BSPFA}$. 

Our Theorem \ref{main} can also be construed as a negative result 
on iterating stationary preserving forcings. (Such negative results
have also been proven long ago by Shelah.) 

\section{The proof.}

\begin{defn}
Let $f$, $g$ both be functions from $\omega_1$ to $\omega_1$.
We shall write $f<^* g$ iff there is some club $C \subset \omega_1$
such that for all $\nu \in C$, $f(\nu) < g(\nu)$.
\end{defn}

Of course, $<^*$ is a well-founded relation on the set of all $f \colon 
\omega_1
\rightarrow \omega_1$. We shall prove Theorem \ref{main}
by showing that ${\sf BMM}$ gives an infinite $<^*$-descending chain
of such functions unless there are inner models with strong cardinals.

In what follows, if $X$ is a set of ordinals
such that there is no inner model
with a strong cardinal containing $X$ then $K(X)$ denotes the core model
{\em over} $X$ (i.e., with $X$ ``thrown in at the bottom''),
and for ordinals $\xi$, $K(X)||\xi$ denotes $K(X)$ cut off at $\xi$.
The reader who is ignorant of the theory of $K$ may always pretend that
$X^\#$ does not exist, in which case $K(X) = L[X]$ and $K(X)||\xi =
L_\xi[X]$; 
of course, doing so only gives a proof of Theorem \ref{main}
where ``for every
$X \in V$ there is an inner model with a strong cardinal containing
$X$'' is replaced by ``for every
$X \in V$, $X^\#$ exists.''

\begin{defn}
Let $a \subset \omega$ be such that 
there is no inner model with a strong cardinal containing $a$, and assume
that
$\omega_1^{K(a)} =
\omega_1^V$. 
Suppose in fact that there are (unique)
$A \subset \omega_1$ and $(a_\nu \colon \nu<\omega_1)$ such that for all
$\nu<\omega_1$, $a_\nu$ is the
$K(A \cap \nu)$-least subset of $\omega$ which is
almost disjoint from each member of $\{ a_{\bar \nu} \colon {\bar \nu}
< \nu \}$, and
$\nu \in A$ iff $a_\nu \cap a$ is finite.

Then we shall denote by $f_a$ the following function:
${\rm dom}(f_a) = \omega_1$, and
for $\nu<\omega_1$, $f_a(\nu) =$ the least $\beta<\omega_1$ such that
$$K(A \cap \nu)||(\beta+1) \models \nu {\rm \ is \ countable.}$$

In this situation, we shall say that $f_a$ exists (or, that
$f_a$ is well-defined). If there are no
$A$, $(a_\nu \colon \nu<\omega_1)$ as above then $f_a$ does not exist.
\end{defn}

Our key lemma is the following.

\begin{lemma}\label{key}
Let $a \subset \omega$ be such that 
there is no inner model with a strong cardinal containing $a$, and assume
that
$f_a$ is well-defined.
There is then a stationary preserving set-generic
extension of $V$ in which there is some $b \subset \omega$, 
$a <_T b$, such that $f_b$ is well-defined and $f_b <^* f_a$. 
\end{lemma}

{\sc Proof} of Theorem \ref{main} from Lemma \ref{key}. Suppose that 
${\sf BMM}$ holds but that for some $X \in V$, there is no inner model with a
strong cardinal containing $X$. We have shown in \cite{K} that there is 
then a stationary preserving set-generic extension of $V$
in which there is some $a \subset \omega$ with $X \in H_{\omega_2} =
K(a)||\omega_2$ (where $\omega_2$ denotes
the $\omega_2$ of the extension).
In this extension, thus
$$\exists a \ \exists 
{\cal M} \ \exists {\cal M}' \ ( 
{\cal M} \models
{\rm \ ``I \ am \ the \ stack \ of \ }
a{\rm -mice \ projecting \ to \ } \omega {\rm ," \ }$$
$${\cal M} \cap {\rm OR} = \omega_1 {\rm , \ } {\cal M}' {\rm \ is
\ transitive \ and \ contains \ all \ sets \ }$$ $${\rm which \
are \ boldface \ definable
\ over \ } {\cal M} {\rm , \ and \ } {\cal M}' \models f_a 
{\rm \ exists \ }).$$
By ${\sf BMM}$, the displayed
statement holds in $V$. If $a_0$, ${\cal M}$, ${\cal M}' \in V$ witness this
then by ${\cal M} \cap {\rm OR} = \omega_1$ and absoluteness,
${\cal M} = K(a_0)||\omega_1$. Moreover, ${\cal M}' \models f_a$ exists
will imply that $f_a$ really exists.

Now let ${\cal C}$ denote the cone of all reals $b$ above $a_0$ in the Turing 
degrees for which $f_b$ exists, 
i.e., ${\cal C} = \{ b \subset \omega : a_0 \leq_T b \wedge f_b$ exists $\}$.
Let $a \in {\cal C}$. By Lemma \ref{key}, there
is a stationary preserving set-generic extension of $V$
in which there is some $b \subset \omega$ with $a <_T b$ and
$f_b <^* f_a$. In this extension, thus
$$\exists b \ \exists {\cal M} \ \exists {\cal M}' \ ( a <_T b
{\rm , \ } {\cal M} \models
{\rm \ ``I \ am \ the \ stack \ of \ }
b{\rm -mice \ projecting \ to \ } \omega {\rm ," \ }$$
$${\cal M} \cap {\rm OR} = \omega_1 {\rm , \ } {\cal M}' {\rm \ is
\ transitive \ and \ contains \ all \ sets \ }$$ $${\rm which \
are \ boldface \ definable
\ over \ } {\cal M} {\rm , \ and \ } {\cal M}' \models f_b <^* f_a).$$
By ${\sf BMM}$, the displayed
statement holds in $V$.
If $b$, ${\cal M}$, ${\cal M}' \in V$ witness this
then ${\cal M} = K(b)||\omega_1$ and ${\cal M}' \models f_b <^* f_a$.
But then $f_b <^* f_a$ really holds true.

But this shows that $<^*$ is not well-founded (in a strong sense:
for each $a \in {\cal C}$, $<^* \upharpoonright \{ f_b 
\colon  
a <_T b \wedge f_b <^* f_a \}$ is ill-founded).
Contradiction! \hfill $\square$ (Theorem \ref{main})
 
\bigskip
{\sc Proof} of Lemma \ref{key}.
Fix $a \subset \omega$ as in the statement of Lemma \ref{key}.
%
Let us 
fix $A \subset \omega_1$, the subset of $\omega_1$ obtained by
``decoding'' $a$. W.l.o.g., $H_{\omega_2}
= K(A)||\omega_2$ (cf.~\cite{K}). 

Let ${\mathbb P} \in V$ be the set of all $(f,c)$ such that there 
is some $\nu < \omega_1$ with: 

\bigskip
\noindent $\bullet \ \ \ $ $f \colon \nu \rightarrow 2$,

\noindent $\bullet \ \ \ $ $c \subset \nu+1$ is closed,

\noindent $\bullet \ \ \ $ for all ${\bar \nu} \leq \nu$, 
$$K(A \cap {\bar \nu},f \upharpoonright {\bar \nu}) \models
{\bar \nu} {\rm \ is \ countable,}$$

\noindent $\bullet \ \ \ $ for all ${\bar \nu} \in c$, 
$$K(A \cap {\bar \nu},f \upharpoonright {\bar \nu})||f_a({\bar \nu}) \models
{\bar \nu} {\rm \ is \ countable.}$$

\bigskip
If $p = (f,c) \in {\mathbb P}$ then we shall write $p^\ell$ for
$f$ and $p^\rho$ for $c$.
A condition $q$ is stronger than $p$ iff $q^\ell \upharpoonright {\rm 
dom}(p^\ell) = p^\ell$ and $q^\rho \cap ({\rm max}(p^\rho)+1) = p^\rho$.

The following is easy to verify.

\bigskip
{\bf Claim 1.} (Extendability)
Let $p \in {\mathbb P}$. If $\nu < \omega_1$ then there is some
$q \leq p$ such that ${\rm dom}(q^\ell) \geq \nu$. Also,
if $\nu < \omega_1$ then there is some $q \leq p$ such that 
$q^\rho \setminus \nu \not= \emptyset$.

\bigskip
Whereas it can be shown that ${\mathbb P}$ is not semi-proper in
general,\footnote{Hint: Otherwise $\forall a \ \exists b
\ f_b <^* f_a$ would hold in the model of \cite{gosh}
if this model is constructed 
by forcing over $L$.} the following does hold true.

\bigskip
{\bf Claim 2.} ${\mathbb P}$ is stationary preserving.

\bigskip
{\sc Proof} of the Claim. Suppose that $p \forces {\dot C} \subset {\check
\omega_1}$ is club, and let $S \subset \omega_1$ be stationary.
We aim to find some $q \leq p$ with $q \forces {\dot C} \cap {\check S}
\not= \emptyset$.

Let $n_0 \in \omega$ be large enough. Let us first pick
$$\pi \colon {\bar K}^* \rightarrow K(A)||\omega_2$$ such that
${\bar K}^*$ is countable and transitive, ${\rm crit}(\pi) \in S$,
and $\{ a, {\mathbb P}, p, {\dot C} \} \subset {\rm ran}(\pi)$.
Set $\nu = {\rm crit}(\pi)$, ${\bar {\mathbb P}} = \pi^{-1}({\mathbb 
P})$, and ${\bar {\dot C}} = \pi^{-1}({\dot C})$.
Working in ${\bar K}^*$ (a model of ${\sf ZFC^-}$), we may pick some 
$${\bar K} \prec_{\Sigma_{n_0}} {\bar K}^*$$
such that ${\bar K} \triangleleft {\bar K}^*$ (i.e.,
the former is a strict initial segment of the latter),
$\rho_{n_0}({\bar K}) = \nu$,
and $\{ a, {\bar {\mathbb P}}, p, {\bar {\dot C}} \} \subset {\bar K}$.
(We may for instance let ${\bar K}$ be the $\Sigma_{n_0}$ hull
of $\nu \cup \{ a, {\bar {\mathbb P}}, p, {\bar {\dot C}} \}$ formed inside
${\bar K}^*$.) We'll have ${\bar K}^* \triangleleft K(A \cap \nu)$.

Set $\beta = {\bar K}^* \cap {\rm OR}$. 

\bigskip 
{\bf Subclaim.} $\beta \leq 
f_a(\nu)$.

\bigskip
{\sc Proof} of the Subclaim. Of course, $\nu$ is uncountable in
${\bar K}^*$, and thus $\nu$ is uncountable in $K(a)^{{\bar K}^*}$.
But a straightforward coiteration argument yields
$K(a)^{{\bar K}^*} \triangleleft K(a)$, i.e., 
$K(a)^{{\bar K}^*} = K(a)||\beta$. Therefore, $\nu$ is uncountable
in $K(a)||\beta$ and hence $\beta \leq f_a(\nu)$. \hfill $\square$ (Subclaim)

\bigskip
We shall now imitate an argument of \cite{shst}.
Let $(E_i \colon i<\nu) \in {\bar K}^*$ be an enumeration 
of all the sets which are club in $\nu$ and which exist in
${\bar K}$, and let $E$ be the diagonal intersection
of $(E_i \colon i<\nu)$. Notice that
$E \setminus E_i$ is bounded in $\nu$ whenever $i<\nu$.
Let us pick an external sequence $(\nu_n \colon n<\omega)$
of ordinals smaller than $\nu$ which is cofinal in $\nu$.
Also, let $\{ D_n \colon n<\omega \}$ be the set of all
sets in ${\bar K}$ which are open dense in ${\bar {\mathbb P}}$.

We now construct a sequence $(p_n \colon n<\omega)$ of conditions
such that $p_0=p$, $p_{n+1} \leq p_n$, and $p_{n+1} \in D_n$
for $n<\omega$. Simultaneously, we'll construct a
sequence $(\delta_n \colon n<\omega)$ of ordinals.

Suppose that $p_n$ is given. Notice that, setting
$\gamma = {\rm dom}(p_n^\ell)$,
$\gamma 
< \nu$ (as $p_n \in {\bar K}$). 
Work inside ${\bar K}$ for a second. Using Claim 1, for all
$\delta$ with $\gamma \leq \delta < \nu$ we may easily pick some 
$p^\delta \leq p_n$ such that: $p^\delta \in D_n$, 
${\rm dom}((p^\delta)^\ell) > {\rm max}(\{ \nu_n , \delta \})$,
and for all limit ordinals $\lambda$ with $\gamma \leq \lambda
\leq \delta$, $(p^\delta)^\ell(\lambda)=1$ iff $\lambda = \delta$.
There is some ${\bar E} \in {\cal P}(\nu) \cap {\bar K}$
club in $\nu$ such that for
any $\eta \in {\bar E}$, $\delta < \eta \Rightarrow {\rm dom}((p^\delta)^\ell)
< \eta$.

Now working inside ${\bar K}^*$, we may pick some 
$\delta \in E$ such that $E \setminus {\bar E} \subset \delta$.
Let us set $p_{n+1} = p^\delta$, and put $\delta_n = \delta$.
Of course, $p_{n+1} \leq p_n$ and $p_{n+1} \in D_n$.
Moreover, ${\rm dom}((p_{n+1})^\ell) < {\rm min}(E \setminus (\delta_n+1))$,
so that for all limit ordinals
$\lambda \in E \cap ({\rm dom}((p_{n+1})^\ell) \setminus
{\rm dom}((p_{n})^\ell))$ we have that
$(p_{n+1})^\ell = 1$ iff $\lambda = \delta_n$.

Now let us define an object $q = (q^\ell,q^\rho)$ as follows.
We set $q^\ell = \bigcup_{n<\omega}(p_n)^\ell$ and
$q^\rho = \bigcup_{n<\omega}(p_n)^\rho \cup \{ \nu \}$. 

Let us verify that $q \in {\mathbb P}$. Well, 
by Claim 1, ${\rm dom}(q^\ell) = \nu$ and $q^\rho \cap \nu$ is unbounded
in $\nu$. Hence to prove that $q \in {\mathbb P}$ boils down to having
to
show that $$K(A \cap \nu,q^\ell) \models
\nu {\rm \ is \ countable.}$$
However, by the construction of the $p_n$'s we have that
$$\{ \lambda \in E \cap ({\rm dom}(q^\ell) \setminus
{\rm dom}(p^\ell)) \colon \lambda {\rm \ is \ a \ limit \ ordinal \
and \ } q^\ell(\lambda) = 1 \} = \{ \delta_n \colon n<\omega \} {\rm , }$$
which is cofinal in $\nu$. But $E \in {\bar K}^* = K(A \cap \nu)||\beta$, 
and therefore
$E \in K(A \cap \nu)||f_a(\nu)$ by the above Subclaim.
Therefore, $\{ \delta_n \colon n<\omega \} \in K(A \cap \nu)||f_a(\nu)$ 
witnesses that $\nu$ is countable in $K(A \cap \nu)||f_a(\nu)$, as desired.

It is now easy to see that $q \forces {\check \nu} \in {\dot C} \cap {\check S}$.
\hfill $\square$ (Claim 2)

\bigskip
The rest is smooth. Let us confuse $V^{{\mathbb P}}$ with 
a generic extension of $V$.
Because forcing with ${\mathbb P}$ does not collapse $\omega_1$,
it adds a pair $B$, $C$ such that $B \subset \omega_1$, $C$ is a 
club subset of $\omega_1$, for all $\nu < \omega_1$,
$$K(A \cap \nu,B \cap \nu) \models \nu {\rm \ is \ countable, }$$ and
for all $\nu \in C$,
$$K(A \cap \nu,B \cap \nu)||f_a(\nu) \models \nu {\rm \ is \ countable.}$$ 
Let us fix 
such a pair $(B,C)$, and let us write $D = A \oplus B$.
Let us code $D$ down to a real in the usual way (cf.~\cite{K}).
In order to do this, let us write $(a_\beta \colon \beta<\omega_1)$
for that sequence of subsets of $\omega$
such that for each $\beta<\omega_1$, $a_\beta$ is the 
$K(D \cap \beta)$-{\em least}
subset of $\omega$ which is almost disjoint from every member
of $\{ a_{\bar \beta} \colon {\bar \beta}<\beta \}$.

Specifically, let ${\mathbb A}$ consist of all pairs $(l(p),r(p))$, where
$l(p) \colon n \rightarrow 2$ for some $n<\omega$ and $r(p) \subset \omega_1$
is finite. A condition $q$ is stronger than $p$ iff $l(q)$ extends $l(p)$,
$r(p)$ is a subset of $r(q)$, and
for all $\beta \in r(q)$, if
$\beta \in D$ then $$\{ n \in {\rm dom}(l(q))\setminus {\rm dom}(l(p)) \colon
l(q)(n) = 1 \} \cap a_\beta = \emptyset.$$
The forcing ${\mathbb A}$ has the c.c.c., and forcing with ${\mathbb A}$
adds a real $b$ such that for all $\beta<\omega_1$,
$$\beta \in D {\rm \ iff \ } b \cap a_\beta {\rm \ is \ finite.}$$

Let us now look at $f_b$. 
Let $C' = \{ \nu \in C \colon K(b)||\nu
\prec_{\Sigma_\omega} K(b)||
\omega_2 \}$.
Of course, $C'$ is club in $\omega_1$.
The proof of the following claim will therefore finish
the proof of Theorem \ref{main}, as $V^{{\mathbb P}
* {\mathbb A}}$ will be an extension as desired.

\bigskip
{\bf Claim 3.} For all $\nu \in C'$, $f_b(\nu) < f_a(\nu)$. 

\bigskip
{\sc Proof} of Claim 3. 
By the choice of $A$, $\nu$ is uncountable in $K(A \cap \nu)||f_a(\nu)$.
However, $\nu$ is countable in $K(D \cap \nu)||f_a(\nu)$.
But $D$ is exactly the subset of $\omega_1$ obtained by ``decoding''
$b$.
Therefore,
we must have $f_b(\nu) < f_a(\nu)$.  
\hfill $\square$ (Claim 3)

\hfill $\square$ (Lemma \ref{key})

\section{A conjecture.}

We do not know how to prove the following.

\bigskip
{\bf Conjecture.} If ${\sf BMM}$ holds then there is an inner model with a 
Woodin cardinal.

\bigskip
In fact, we do not even know how to get $0^\P$ from ${\rm BMM}$.
This is related to the problem that we do not know 
how to get $0^\P$ from the assumption that the theory of $L({\mathbb R})$
is absolute for stationary preserving forcings (cf.~\cite{K}).

\end{document}